\documentclass[a4paper,11pt]{article}
\usepackage{graphicx}
\usepackage{amsmath,amsthm}
\usepackage{amsfonts,amssymb,stmaryrd}
\usepackage{amscd}
\usepackage[usenames,dvipsnames,svgnames,table]{xcolor}
\usepackage[cp1251]{inputenc}
\usepackage[matrix,arrow,curve]{xy}
\usepackage{latexsym}
\usepackage{fullpage}
\usepackage[breaklinks, naturalnames]{hyperref}

\newtheorem{lemma}{Lemma}

\newtheorem{theorem}[lemma]{Theorem}
\newtheorem{coro}[lemma]{Corollary}

\newtheorem{thm}{Theorem}
\newtheorem{Def}[thm]{Definition}

\newtheorem{ex}{Example}

\newcommand{\Hom}{{\rm Hom}}
\newcommand{\uHom}{{\rm hom}}
\newcommand{\End}{{\rm End}}

\newcommand{\kk}{{\bf k}}

\newcommand{\Ext}{{\rm Ext}}
\newcommand{\Tor}{{\rm Tor}}
\newcommand{\uExt}{{\rm ext}}
\newcommand{\Mod}{{\rm Mod}}

\newcommand{\bZ}{{\mathbb Z}}

\newcommand{\SSS}{{\mathcal S}}
\newcommand{\RR}{{\mathcal R}}

\newcommand{\la}{\langle}
\newcommand{\ra}{\rangle}
\newcommand{\N}{\mathbb{N}}

\newcommand{\ot}{\otimes}

\newcommand{\tv}{\tilde v}

\renewcommand{\le}{\leqslant}
\renewcommand{\ge}{\geqslant}

\newenvironment{Proof}[1][Proof.]{\begin{trivlist}
\item[\hskip \labelsep {\bfseries #1}]}{\flushright
$\Box$\end{trivlist}}

\numberwithin{equation}{section}
\begin{document}
\title{Generating degrees for graded projective resolutions}
\author{E. Marcos, A. Solotar, Y. Volkov
\thanks{\footnotesize The first named author  has been supported by the thematic project of Fapesp
2014/09310-5. The second named author has been partially  supported by projects PIP-CONICET 11220150100483CO and 
UBACyT 20020130100533BA. The first and second authors have been partially supported by project MathAmSud-REPHOMOL.
The third named author has been supported by a post-doc scholarship of  Fapesp 
(Project number:  2014/19521-3) and by Russian Federation Presedent grant (Project number: MK-1378.2017.1).
The second named author is a	research member of CONICET (Argentina).}}
\date{}
\maketitle
\begin{abstract}
We provide a framework connecting several well known theories related to the linearity of graded modules over graded algebras. In the first part, we pay a particular attention to the tensor products of graded bimodules over graded algebras.
Finally, we provide a tool to evaluate the possible degrees of a module appearing in a graded projective resolution once the generating degrees for the first term of some particular projective resolution are known. 
\end{abstract}

\noindent 2010 MSC: 16S37, 18G10.

\noindent Keywords: Koszul, linear modules, Gr\"obner bases.

\section{Introduction}
 Koszul algebras were introduced by S.~Priddy  in \cite{pr}.
We will apply the notion of a Koszul algebra for algebras presented by quivers with relations. It can be stated as follows. Suppose that $\kk$ is a field, $Q$ is a finite quiver, $I$ is a homogeneous 
ideal of the path algebra $\kk Q$ and 
 $A= \kk Q/I$. 
The algebra $A$ is called {\em Koszul} if the maximal semisimple graded quotient $A_0$ of $A$ has a graded $A$-projective 
resolution $(P_{\bullet}, d_{\bullet})$ such that
for all $i\ge 0$, the $A$-module $P_i$ is generated in degree $i$.
Such a resolution is called a linear resolution and it is minimal whenever it exists, in the sense that $d_i(P_{i+1})\subset P_iA_{>0}$ for all $i\ge 0$.
 
E.~Green and R.~Mart\'\i nez Villa proved in \cite{gm} that the quadratic algebra $A$ is Koszul if and only if its Yoneda algebra, 
$E(A)= \oplus_{i\ge 0} \Ext^i_A(A_0, A_0)$ is generated in degrees 0 and 1, which in turn is equivalent to the Yoneda algebra being isomorphic to the quadratic dual  
$A^{!}$ of $A$. 

Koszulness has been generalized to various settings. Next we describe some of these generalizations.

R.~Berger introduced in \cite{be} the notion of ``nonquadratic Koszul algebra" for algebras of the form $A=T_{\kk}V/I$, 
where $V$ is a finite dimensional $\kk$-vector space and $I$ is a two-sided ideal
generated in degree $s$, for some $s \geq 2$. He required the trivial $A$-module $\kk$ to have a minimal graded projective resolution $(P_{\bullet}, d_{\bullet})$ such that  each $P_i$ is 
generated in degree $\frac{is}{2}$ for $i$ even and $\frac{(i-1)s}{2} + 1$ for $i$ odd.

The authors of \cite{gmmz} considered, under the name of $s$-Koszul algebras, non necessarily quadratic Koszul algebras of 
the form $A=\kk Q/I$,  with $Q$ a finite quiver and $I$ an ideal generated by 
homogeneous elements of degree $s$, connecting this notion with the Yoneda algebra: the algebra $A$ is $s$-Koszul if and only if $E(A)$ is generated in degrees 0, 1 and 2. 
Observe that $2$-Koszul algebras are just Koszul algebras. 

Later on, E.~Green and E.~Marcos generalized this notion defining $\delta$-Koszul and $\delta$-determined algebras. See \cite{gm1} for details.

Moreover, E.~Green and E.~Marcos also introduced in \cite{gm2} a family of algebras that they called $2$-$s$-Koszul. 
They proved that these algebras also have the property that their Yoneda algebras are generated
in degrees 0, 1, and 2.

The main objective of the current work is to place all these definitions in a unique framework. We next sketch how we will do this. 

Let $A=\oplus_{i\ge 0} A_i$ be a graded $\kk$-algebra generated in degrees 0 and 1, such that $A_0$ is a finite direct product of fields and $A_1$ is finite dimensional. 
Given a graded $A$-module $X$, we consider a minimal graded projective resolution $(P_{\bullet}(X), d_{\bullet}(X))$ and we take into account 
in which degrees $P_i(X)$ is generated 
for each $i\ge 0$. We are specially interested in what we call $\SSS$-determined case.

In Section \ref{spectralsequence} In Section \ref{spectralsequence} we prove that, given graded $\kk$-algebras $A,B$ and $C$,  a graded 
$A-B$ bimodule $X$ and  a graded $B-C$ bimodule $Y$, if $X$ has a linear minimal $A-B$ -projective
graded resolution, $Y$ has a linear minimal $B-C$ -projective graded resolution, and $\Tor_i^B(X,Y)$ vanishes for $i\ge 1$, then $X\ot_BY$ has a linear minimal $A-C$ -projective graded resolution. 
This is a particular case of  Theorem \ref{main} below. 
 Note that this theorem shows that any graded bimodule over a Koszul algebra which is linear as a right module and flat as a left module is also linear as a bimodule and the tensor product with such a module gives a functor 
from the category of linear graded modules to the category of linear graded modules. Moreover, it recovers and generalizes the fact that the tensor product of two Koszul algebras is Koszul.
The same holds for the $\SSS$-determined (see Definition \ref{Slin}).

 Section \ref{grobner} is devoted to Gr\"obner bases. 
Loosely speaking, we show how one can use them to obtain generating degrees for the $n$-th term of the minimal graded projective resolution of a module if one knows the generating degrees for 
terms of its projective presentation of a special form.

We fix a field $\kk$. All algebras will be $\kk$-algebras and all modules will be right $A$-modules unless otherwise stated. We will simply write $\otimes$ for $\otimes_{\kk}$ and 
$\mathbb{N}_0$ for the set of non negative integer numbers.

We thank the referee for the suggestions and for a careful reading of a previous version of this paper.
 
\section{Tensor products of $\SSS$-determined modules}
\label{spectralsequence}

In this section we will prove Lemma \ref{gr_spec}, which is a graded version of the spectral 
sequences (2) and (3) from \cite[page 345]{CE}.

Let $A, B$ and $C$ be $\kk$-algebras. Let $X$ be an $A-B$-bimodule, $Y$ a $B-C$-bimodule and $Z$ an $A-C$-bimodule. 
Given $a\in A$, we will denote left multiplication by $a$ on $X$ by $L_a\in\End_B(X)$.
We recall that for each $n\in \mathbb{N}$, $\Ext^n_C(Y,Z)$ 
is an $A-B$ bimodule with the structure given by 
\[
aTb:=\big(\Ext^n_C(L_b,Z)\circ \Ext^n_C(Y,L_a)\big)(T), \hbox{ for } T\in \Ext^n_C(Y,Z), a\in A, b\in B.
\]

Suppose now  that $A$ is a $\bZ$-graded algebra and $M$ is a graded $A$-module.
Given $i\in \bZ$, $M[i]$ will denote the $i$-shifted graded $A$-module with underlying $A$-module structure as before, 
whose grading is such that $M[i]_r=M_{i+r}$.
For any graded $A$-module $N$ and any $n\in \N$, we will denote by $\Hom_{Gr A}(M,N)$ the set of degree preserving $A$-module maps from $M$ to $N$
and by  $\Ext^n_{Gr A}(M,N)$ the set of equivalence classes of exact sequences of 
graded $A$-modules  with degree zero morphisms
\[
0\rightarrow N\xrightarrow{f_{n-1}}T_{n-1}\xrightarrow{f_{n-2}}\dots\xrightarrow{f_0}T_0\xrightarrow{f_{-1}}M \rightarrow 0.
\]
Let us consider as usual $\uExt_A^n(M,N):=\oplus_{i\in\mathbb{Z}}\Ext^n_{Gr A}(M,N[i])$, which is a subset of $\Ext_A^n(M,N)$ in a natural way.
Moreover, if $M$ has an $A$-projective resolution with finitely generated  modules, then 
both sets coincide. 

Suppose now that $A$, $B$ and $C$ are $\bZ$-graded algebras, and that the bimodules  $X,Y$ and $Z$ are graded.
For each $n\ge 0$, the $A-B$-bimodule structure on $\Ext^{n}_C(Y,Z)$ induces a graded $A-B$-bimodule structure on $\uExt_C^n(Y,Z)$ whose $i$-th component 
is  $\Ext^n_{Gr C}(Y,Z[i])$. Note also that  $\uExt_C^n(Y,Z[i])\cong \uExt_C^n(Y,Z)[i]$ as graded $A-B$-bimodule, moreover  for any $n\ge 0$, $\Tor_n^B(X,Y)$ is a graded $A-C$-bimodule in a natural way.

The main tool of this section is the following lemma.

\begin{lemma}\label{gr_spec} Let $A$, $B$ and $C$ be $\bZ$-graded algebras, $X$ a graded $A-B$-bimodule, $Y$ a graded $B-C$-bimodule, and $Z$ a graded $A-C$-bimodule.
There are two first quadrant cohomological spectral sequences with second pages
\[
E^{i,j}_2=\Ext^i_{Gr(A^{op}\ot B)}\big(X,\uExt^j_C(Y,Z)\big)\mbox{ and }\tilde E^{i,j}_2=\Ext^i_{Gr(A^{op}\ot C)}\big(\Tor_j^B(X,Y),Z\big)
\]
that converge to the same graded space.
\end{lemma}
\begin{Proof}
 Let
\[
\dots \xrightarrow{d_n(X)} P_n(X)\xrightarrow{d_{n-1}(X)}\dots \xrightarrow{d_0(X)}P_0(X)(\xrightarrow{\mu_X}X)
\]
be a graded $A-B$-projective resolution of $X$ and
\[
(Z\xrightarrow{\iota_Z})I^0(Z)\xrightarrow{d^0(Z)}\dots\xrightarrow{d^{n-1}(Z)}I^n(Z)\xrightarrow{d^n(Z)}\dots
\]
be a graded $A-C$-injective resolution of $Z$. 
Consider two bicomplexes whose $(i,j)$-components are respectively 
\[
\Hom_{Gr(A^{op}\ot B)}\Big(P_i(X),\uHom_C\big(Y, I^j(Z)\big)\Big)
\]
and 
\[
\Hom_{Gr(A^{op}\ot C)}\big(P_j(X)\otimes_BY, I^i(Z)\big).
\]
Since there is an isomorphism of complexes
$$
F_{\bullet}=\Hom_{Gr(A^{op}\ot C)}\big(P_{\bullet}(X)\otimes_B Y, I^{\bullet}(Z)\big)\cong \Hom_{Gr(A^{op}\ot B)}\Big(P_{\bullet}(X),\uHom_C\big(Y, I^{\bullet}(Z)\big)\Big), 
$$
the respective total complexes are isomorphic.
Here, as usually, for two complexes of graded modules $\big(U_{\bullet},d_{U,\bullet}\big)$ and $\big(V^{\bullet},d^{V,\bullet}\big)$ over the algebra $D$ we denote by $\Hom_{Gr D}(U_{\bullet},V^{\bullet})$ the complex
with $\big(\Hom_{Gr D}(U_{\bullet},V^{\bullet})\big)_n=\oplus_{i\in\bZ}\Hom_{Gr D}(U_{i-n},V^{-i})$ and differential $d_{\bullet}$ defined by the equality $d_n(f)=d^{V,-i}f+(-1)^nfd_{U,i-n}$ for $f\in\Hom_{Gr D}(U_{i-n},V^{-i})$.

The first two pages of the spectral sequence $E$ corresponding to the first bicomplex are
\[
E^{i,j}_1=\Hom_{Gr(A^{op}\ot B)}\big(P_i(X),\uExt^j_C(Y,Z)\big)\mbox{ and }E^{i,j}_2=\Ext^i_{Gr(A^{op}\ot B)}\big(X,\uExt^j_C(Y,Z)\big),
\]
while the first two pages of the spectral sequence $\tilde E$ corresponding to the second bicomplex are
\[
\tilde E^{i,j}_1=\Hom_{Gr(A^{op}\ot C)}\big(\Tor_j^B(X,Y), I^i(Z)\big)\mbox{ and }\tilde E^{i,j}_2=\Ext^i_{Gr(A^{op}\ot B)}\big(\Tor_j^B(X,Y),Z\big).
\]
Since both spectral sequences converge to the homology of $F_{\bullet}$, the lemma is proved.
\end{Proof}

From now on any $\bZ$-graded algebra $A$ is assumed to be non negatively graded, that is $A=\oplus_{i\ge 0}A_i$, where $A_0$ is 
isomorphic to a finite product of copies of $\kk$ as an algebra, $\dim_{\kk}A_1<\infty$, and $A$ is generated as an algebra by 
$A_0\oplus A_1$. This is equivalent to say that $A\cong (\kk Q)/I$ where $Q$ is a finite quiver and $I$ is an ideal generated by homogeneous elements of 
degree bigger or equal 2.

\begin{Def}\label{Slin}

{\rm
Let $\SSS=(\SSS_i)_{i\ge 0}$ be a collection of subsets $\SSS_i\subset \bZ$. A graded $A$-module $X$ is called {\it $\SSS$-determined} if it has a 
graded projective resolution $P_{\bullet}(X)$ such that $P_i(X)$ is generated as $A$-module by elements of degrees belonging to
$\SSS_i$, i.e. $P_i(X)=\la\oplus_{j\in \SSS_i}P_i(X)_j\ra_A$ for all $i\ge 0$.  We say that $X$ is {\it $\SSS$-determined up to  degree $r$} if the condition on $P_i(X)$ holds for $0\le i\le r$.}\\

If the set $\SSS_i =\{i\}$, i.e. each $P_i(X)$ is generated in degree $i$, then we say that the resolution is linear.
\end{Def}

Equivalently, a graded $A$-module $X$ is $\SSS$-determined if and only if for any $i\ge 0$ and any graded $A$-module $Y$ with support
not intersecting $\SSS_i$ -- that is, $\oplus_{j\in \SSS_i}Y_j=0$ -- the space $\Ext^i_{Gr A}(X,Y)$ is zero. Analogously, the graded $A$-module $X$ is $\SSS$-determined up to degree $r$ if and only if the last mentioned condition holds for $0\le i\le r$.

 The notion of an $\SSS$-determined module provides a general framework for some well-known situations. We will now exhibit some well known examples of
 $\SSS$-determined  modules. 

\begin{itemize}
 \item  Consider a function $\delta:\bZ_{\ge 0}\rightarrow \bZ$, and define $\SSS_i=\{\delta(i)\}$ for all $i\ge 0$, the $\SSS$-determined 
modules are called {\it $\delta$-determined modules}. If $A_0$ is a $\delta$-determined module over $A$, 
then the graded algebra $A$ is called 
$\delta$-determined.
\item  With the same notations, if moreover the $\Ext$ algebra, $E(A)$, of $A$ is finitely generated, then $A$ is called
{\it $\delta$-Koszul}, see \cite{gm1}. 
In particular, if $\delta$ is the identity, then $\delta$-determined modules are called {\it linear modules} and $\delta$-Koszul 
algebras are exactly {\it Koszul algebras} \cite{pr}. 
\item   Also, given $s\in \mathbb{N}$, let us define $\chi_s:\mathbb{N}_0\rightarrow \bZ$ by
\begin{equation*}
     \chi_s(i)=\left\lbrace
\begin{array}{ll}
     \frac{is}{2} &\text{ if $i$ is even,} \\
     \frac{(i-1)s}{2}+1 &\text{ if $i$ is odd.}\\
\end{array}
\right.
\end{equation*}
The $\chi_s$-linear modules are called {\it $s$-linear modules} and $\chi_s$-Koszul algebras are {\it $s$-Koszul algebras}, see \cite{be}. 
Denoting by $\SSS_i$ the set $\{j\mid j\le\chi_s(i)\}$, $\SSS$-linear modules correspond to {\it $2$-$s$-linear modules}.
If moreover $A_0$ is a $2$-$s$-determined module over $A$, then the graded algebra $A$ is called {\it $2$-$s$-determined}, see \cite{gm2}.
\end{itemize}

Using minimal graded projective resolutions, it is not difficult to see that the $A$-module  $A_0$ 
is $\SSS$-determined if and only if $A$ is an $\SSS$-determined module over $A^{op}\ot A$. 
This fact follows, for example, from \cite[Theorem 2]{Skol}.

Given two collections $\SSS=(\SSS_i)_{i\ge 0}$ and $\RR=(\RR_i)_{i\ge 0}$ of subsets of $\bZ$ we define the collection 
$\SSS\ot \RR=\big((\SSS\ot \RR)_i\big)_{i\ge 0}$ by 
\[
(\SSS\ot \RR)_i=\bigcup\limits_{\scriptsize\begin{array}{c}j+k=i\\j,k\ge 0\end{array}}\{n+m\mid n\in \SSS_j, m\in \RR_k\}.
\]
Lemma \ref{gr_spec} allows us to prove the following theorem, which generalizes some well known results about Koszul algebras 
and Koszul modules  concerning tensor products.

\begin{thm}\label{main} Let $\SSS=(\SSS_i)_{i\ge 0}$  and $\RR=(\RR_i)_{i\ge 0}$ be two collections of subsets of $\bZ$. 
Let $A$, $B$ and $C$ be $\bZ$-graded algebras,  and finally let $X$ be a graded $A-B$-bimodule which is $\SSS$-determined as bimodule and $Y$ be a graded $B-C$-bimodule which is $\RR$-determined
as $C$-module. 
If $\Tor_i^B(X,Y)=0$ for $1\le i\le r-1$, then $X\otimes_B Y$ is an $\SSS\ot \RR$-determined until $r$-th degree $A-C$-bimodule. In particular, if $\Tor_i^B(X,Y)=0$ for all $i\ge 1$, then $X\otimes_B Y$ is an $\SSS\ot \RR$-determined $A-C$-bimodule.
\end{thm}
\begin{Proof} Let us fix $n$ and $r$ such that $0\le n\le r$. For any graded $A-C$-bimodule $Z$ such that $\oplus_{m\in (\SSS\ot \RR)_n}Z_m=0$,
we will prove that 
$\Ext^n_{Gr(A^{op}\ot C)}(X\ot_BY,Z)=0$. 
By Lemma \ref{gr_spec} there are spectral sequences
\[
E^2_{i,j}=\Ext^i_{Gr(A^{op}\ot B)}\big(X,\uExt^j_C(Y,Z)\big)\mbox{ and }\tilde E^2_{i,j}=\Ext^i_{Gr(A^{op}\ot C)}\big(\Tor_j^B(X,Y),Z\big)
\]
that converge to the same graded space $T_{\bullet}$. It follows easily from the condition on $\Tor_*^B(X,Y)$ that $T_{n}=\Ext^n_{Gr(A^{op}\ot C)}(X\otimes_BY,Z)$ if $n<r$ and that $T_{r}=\Ext^r_{Gr(A^{op}\ot C)}(X\otimes_BY,Z)\oplus V$ for some $V\subset \Hom_{Gr(A^{op}\ot C)}\big(\Tor_r^B(X,Y),Z\big)$.

Thus, it is enough to prove that $E^2_{i,j}=0$ for all integers $i,j\ge 0$ such that $i+j=n$. 
Let us fix such $i$ and $j$. If $k\in \SSS_i$, then for any $l\in \RR_j$ it is clear that $k+l\in (\SSS\ot \RR)_n$ and so 
$Z[k]_l=Z_{ k+l}=0$.
Since $Y$ is an $\RR$-linear $C$-module, we know that $\uExt^j_C(Y,Z)_k=\Ext^j_{Gr C}(Y,Z[k])=0$ for any $k\in \SSS_i$;  from this, 
since $X$ is an $\SSS$-linear $A^{op}\ot B$-module, $E^2_{i,j}=\Ext^i_{Gr(A^{op}\ot B)}\big(X,\uExt^j_C(Y,Z)\big)=0$.
We have proven that for any $0\le n\le r$ and any graded $A-C$-bimodule $Z$ such that $\oplus_{m\in (\SSS\ot \RR)_n}Z_m=0$ one has $\Ext^n_{Gr(A^{op}\ot C)}(X\ot_BY,Z)=0$. 
Consequently, $X\ot_BY$ is an $\SSS\ot \RR$-determined until $r$-th degree $A-C$-bimodule.
\end{Proof}

\begin{ex} Let $A$ be the $\kk$-algebra with generators $x$ and $y$ subject to the relations $xy=yx=0$ and $x^3=y^3$. Let $X=A/\langle x\rangle$ and $Y=A/\langle y\rangle$.
Note that $X$ is a graded $A$-module and $Y$ is a graded $A$-bimodule in a natural way. We will show the the conclusion of the Theorem \ref{main} to the tensor product $X\otimes_AY$ of $\kk$-$A$-bimodule $X$ and $A$-$A$-bimodule $Y$, does not hold.
One can easily see that there are short exact sequences 
\[ Y[1]\hookrightarrow A\twoheadrightarrow X\]
and 
\[X[1]\hookrightarrow A\twoheadrightarrow Y\]
of graded right $A$-modules. It follows, from these two short exact sequences that, $X$ and $Y$ are linear as right $A$-modules. On the other hand, 
$X\otimes_AY\cong A/\langle x,y\rangle$ is the unique simple $A$-module whose 
minimal projective resolution $P_{\bullet}$ is not linear at $P_2$. This example shows that Theorem \ref{main} is not valid without the vanishing condition on $\Tor^B_*(X,Y)$.
\end{ex}

\begin{coro} Let $X$ be a graded $A-B$-bimodule. If $A_0$ is an $\SSS$-determined right
$A$-module and $X$ is an $\RR$-determined $B$-module, 
then $X$ is an $\SSS\ot \RR$-determined $A-B$-bimodule. In particular,
if $A$ is $2$-$s$-determined and $X$ is a $2$-$s$-determined $B$-module, then $X$ is a $2$-$s$-determined $A-B$-bimodule. 
When $s =2$ we obtain that  if $A$ is Koszul and $X$ is a linear $B$-module, then $X$ is a linear $A-B$-bimodule.
\end{coro}
\begin{Proof} Since $A_0$ is an $\SSS$-determined right
 $A$-module, $A$ is an $\SSS$-determined module over $A^{op}\ot A$.
Since $A$ is flat as a right $A$-module, the result follows from Theorem \ref{main}, since $A\ot_AX\cong X$.
\end{Proof}

It is was proved in \cite{bf} and \cite {gm} that if $A$ and $B$ are Koszul algebras, then $A\ot B$ is a Koszul algebra too. In the next corollary we give a very short and easy proof of a generalization of this fact.
Note that this generalization follows from \cite[Chapter 3, Proposition 1.1]{PP} for algebras $A$, $B$ such that $A_0=B_0=\kk$. 

\begin{coro} If $X$ is an $\SSS$-linear $A$-module and $Y$ is an $\RR$-linear $B$-module, then $X\ot Y$ is an 
$\SSS\ot \RR$-linear module over $A\ot B$. In particular, if $A$ and $B$ are
$2$-$s$-determined, then $A\ot B$ is $2$-$s$-determined. In particular, if $A$ and $B$ are Koszul, then $A\ot B$ is Koszul.
\end{coro}
\begin{Proof} It follows from Theorem \ref{main} since any $\kk$-module is flat. The second part follows from the fact that $(A\ot B)_0=A_0\ot B_0$.
\end{Proof}

It is interesting to mention also the following special case of Theorem \ref{main}.

\begin{coro} Let $X$ be a graded $A-B$-bimodule that is flat as a left $A$-module. If $X$ is $2$-$s$-determined as $B$-module, then the functor 
$-\ot_AX:\Mod A\rightarrow\Mod B$ induces a functor from the category of $2$-$s$-determined $A$-modules to the category of 
$2$-$s$-determined $B$-modules. In particular, if $X$ is a linear $B$-module, then $-\ot_AX$ induces  a functor from the category of 
linear $A$-modules to the category of linear $B$-modules
\end{coro}

\section{Using Gr\"obner bases}
\label{grobner}

In this section we will use Gr\"obner bases techniques to study graded projective resolutions of a graded module $X$ over an algebra $A=\kk Q/I$,
where $Q$ is a finite quiver and $I$ is a homogeneous ideal contained in $(\kk Q_{>0})^2$. Our aim is to estimate   the degrees of the modules appearing in the minimal projective 
resolution of $X$ using Gr\"obner basis of $I$ and a particular graded projective presentation of $X$.

We fix a set of paths $S\subset Q_{\ge 2}$. Next we introduce some notation.
Given two paths $p$ and $q$ in $Q$, we write $p\mid q$ if there are paths $u$ and $v$ in $Q$ -- possibly of length $0$ -- such that $q=upv$. If $q=pv$ (resp. $q=up$) we say 
that $p$ divides $q$ on the left (resp. right) and we write $p\mid_l q$ (respectively $p\mid_rq$). We say that $S$ is {\it reduced} if for any $q\in S$ there is no $p \in S$,  $p\neq q$ such that 
$p\mid q$. Let us write $len(q)=n$ if $q\in Q_n$.

We next define, for $n \in \mathbb{N}$, the notion of $n$-{\em overlap} for $S$. These elements will provide a  minimal set of generators for each projective 
module of the minimal projective resolution 
of $A_0$ as $A$-module. Note that $n$-overlaps are called $n$-chains by D.~Anick  in \cite{an} and $n$-ambiguities by S.~Chouhy and A.~Solotar in \cite{cs}.

Given a quiver $Q$ we also denote by $Q$ the set of paths in $Q$, the context makes it clear what we mean.

\begin{Def}{\rm
Let $Q$ be a quiver and $S\in Q_{\ge 2}$ be a set of paths. We say that $p\in Q$ is an {\it $S$-path} if there exists $s\in S$ such that $s\mid_r p$. 
We denote by $Q_S$ the set of $S$-paths in $Q$.
Suppose that $p=qu$ for some $u,q\in Q$. We say that $q$ {\it $S$-vanishes} $p$ if there is no $s\in S$ dividing $u$.
We say that $q$ {\it almost $S$-vanishes} $p$ if  $q$ does not $S$-vanish $p$ and, for any presentation $u=u_1u_2$ with $u_2\in Q_{>0}$, $q$ $S$-vanishes $qu_1$. We write $q\mid_S p$ if $q$ $S$-vanishes $p$ and $q\mid_S^ap$ if $q$ almost $S$-vanishes $p$ (note that the relations $\mid_S$ and $\mid_S^a$ are not transitive). If $q\mid_S^ap$, then we automatically have $p\in Q_S$.

We next define the set of {\it $n$-overlaps} $O_n(S)\subset Q$ and the set of {\it $n$-quasioverlaps} $QO_n(S)\subset Q\times Q_{>0}$ inductively on $n$.
\begin{itemize}
\item For $n=0$, we define $O_0(S)=Q_1$ and $QO_0(S)=\{(w,v)\mid w\in Q_0,v\in Q_{>0}, vw =v\}$.
\item For $n=1$, we define $O_1(S)=S$ and $QO_1(S)=\{(w,v)\mid  w,v\in Q_{>0},vw\in S\}$.
\item For $n>1$, we define 
$$O_n(S)=\{w\mid \exists w_1\in O_{n-1}(S),w_2\in O_{n-2}(S)\mbox{ such that }w_1\mid_Sw, w_2\mid^a_Sw\}$$
and
\begin{multline*}
QO_n(S)\\=\{(w,v)\mid \exists (w_1,v)\in QO_{n-1}(S),(w_2,v)\in QO_{n-2}(S)\mbox{ such that }w_1\mid_Sw, w_2\mid^a_Sw\}.
\end{multline*}
\end{itemize}
}
\end{Def}

\begin{lemma}
Suppose that $S$ is reduced and $n\ge 1$. If $w\in O_n(S)$, then there is unique $w'$ such that $w'\mid_l w$ and $w'\in O_{n-1}(S)$. If $(w,v)\in QO_n(S)$, then there is unique $w'$ such that $w'\mid_l w$ and  $(w',v)\in QO_{n-1}(S)$.
\end{lemma}
\begin{Proof}
We will only prove the assertion about $n$-overlaps since the proof for $n$-quasioverlaps is similar. The existence of $w'$ is a direct consequence of the definition of an $n$-overlap. 
Thus the only thing to prove is uniqueness.
We proceed by induction on $n$. For $n=1$, the assertion is obvious. For $n=2$, the assertion follows directly from the fact that $S$ is reduced.  Suppose now that $n>2$ and we have already proven 
the statement for $n-1$ and $n-2$. Suppose that there are two different paths $w',w''\in O_{n-1}(S)$ such that $w'\mid_lw$ and $w''\mid_lw$. Without loss of generality we may assume that 
$w''=w'u$ for some $u\in Q_{>0}$.
By the definition of $O_{n-1}(S)$, there exist $w'_1\in O_{n-2}(S)$ and $w'_2\in O_{n-3}(S)$ such that $w'_1\mid_Sw'$ and $w'_2\mid_S^aw'$.
The inductive hypothesis assures that $w'_1$ and $w'_2$ are unique elements of $O_{n-2}(S)$ and $O_{n-3}(S)$ respectively such that $w'_2\mid_lw'_1\mid_l w''$. 
Then $w'_2\mid_S^aw''$. Since $u\in Q_{>0}$, we have $w_2'\mid_S w'$. Since this is a contradiction, the proof is complete.
\end{Proof}

Given $w\in O_{n}(S)$, and $i$ such that $0\le i\le n$, $\rho_i(w)$ will denote the unique element of $O_{i}(S)$ that divides $w$ on the left. Analogously, 
for $(w,v)\in QO_{n}(S)$, $0\le i\le n$, $\rho_i^v(w)$ will denote the unique element such that $(\rho_i^v(w),v)\in QO_{i}(S)$ and $\rho_i^v(w)\mid_lw$. The next two lemmas give alternative 
definitions for $n$-overlaps and $n$-quasioverlaps.

\begin{lemma}\label{partition_o}
Consider a reduced set of paths $S$ and an integer $n$ such that $n\ge 1$. Given a path $w\in Q$, $w\in O_n(S)$ if and only if it can be represented in the 
form $w=v_0u_1v_1u_2\dots u_{n-1}v_{n-1}u_nv_n$, where $u_1,\dots, u_n\in Q$, $v_1,\dots, v_{n-1} \in Q_{>0}$, and $v_0,v_n \in Q_0$ are such that
\begin{enumerate}
\item  $v_{i-1}u_iv_i\in S$ for $1\le i\le n$,
\item for all $1\le i \le n-1$, there are no $u,v \in Q_{>0}$ such that  
$vv_iu\in S$, $u\mid_l u_{i+1}v_{i+1}$ and $v\mid_r u_i$,
\item $len(u_1)>0$ if $n\le 2$.
\end{enumerate}
\end{lemma}
\begin{Proof} Suppose that $w\in O_n(S)$. Given $i$, $1\le i\le n$, there exists a unique $w_i\in S$ such that $w_i\mid_r\rho_i(w)$. It follows from the definition of $O_i(S)$ that, 
for $2\le i\le n$, there exist $v_{i-1},u_i'\in Q_{>0}$ such that $\rho_i(w)=\rho_{i-1}(w)u_i'$ and $w_i=v_{i-1}u_i'$. We also define $u_1'=\rho_1(w)$. Using again the definition of $i$-overlap, 
we get $v_i\mid_ru_i'$ for $1\le i\le n-1$. It remains to define $u_i$ from the equality $u_i'=u_iv_i$ for $1\le i\le n-1$ and $u_n=u_n'$. It is clear that $v_0$ and $v_n$ are simply the ending 
and the starting vertices of $w$.

Now, if $u_1,\dots, u_n\in Q$, $v_1,\dots, v_{n-1} \in Q_{>0}$, and $v_0,v_n \in Q_0$ satisfy all the required conditions, then the induction on $1\le i\le n$ shows that $u_1v_1\dots u_iv_i\in O_i(S)$.
\end{Proof}

\begin{lemma}\label{partition_qo}
Let $S$ be a reduced set of paths and let $n$ be an integer, $n\ge 1$. Given paths $w\in Q$ and $v_0\in Q_{>0}$, the element$(w,v_0)\in QO_n(S)$ if and only if $w$ can be represented in the 
form $w=u_1v_1u_2\dots u_{n-1}v_{n-1}u_nv_n$, where $u_1,\dots, u_n\in Q$, $v_1,\dots, v_{n-1} \in Q_{>0}$, and $v_n \in Q_0$ are such that
\begin{enumerate}
\item  $v_{i-1}u_iv_i\in S$ for $1\le i\le n$,
\item for all $1\le i \le n-1$, there are no $u,v \in Q_{>0}$ such that  
$vv_iu\in S$, $u\mid_l u_{i+1}v_{i+1}$ and $v\mid_r u_i$,
\item $len(u_1)>0$ if $n=1$.
\end{enumerate}
\end{lemma}
\begin{Proof} The proof is analogous to the proof of Lemma \ref{partition_o} and so it is left to a reader.
\end{Proof}

\begin{ex}
 Let $Q$ be the quiver with $Q_0=\{e\}$ and $Q_1=\{x,y\}$. Fix $S=\{x^2y^3, x^3\}$. The element $(w_0,v_0)=(xxxyyy,xx)$ is a $3$-quasioverlap with $\rho_2^{v_0}(w_0)=xxx$ and $\rho_1^{v_0}(w_0)=x$. 
The elements $u_1=e$, $v_1=u_2=v_2=x$, and $u_3=yyy$ provide the partition of Lemma \ref{partition_qo}.
At the same time, the element  $w=v_0w_0=xxxxxyyy$ is a $3$-overlap with $\rho_2(w)=xxxx$ and $\rho_1(w)=xxx\in S$. In this case the paths $u_1=x$, $v_1=xx$, $u_2=e$, $v_2=x$, and $u_3=xyyy$ 
provide the partition of Lemma \ref{partition_o}. Note that, though $v_0w_0$ is a $3$-overlap, we have $\rho_2(v_0w_0)\not=v_0\rho_2^{v_0}(w_0)$ and this fact causes differences in the partitions of $v_0w_0$ and $(w_0,v_0)$.
\end{ex}

\begin{ex} Let us take $Q$ as in the previous example and $S=\{x^3,xy^2\}$. Consider $(w_0,v_0)=(xxxyy,xx)$. It is a $3$-quasioverlap with $\rho_2^{v_0}(w_0)=xxx$ and $\rho_1^{v_0}(w_0)=x$. 
The elements $u_1=e$, $v_1=u_2=v_2=x$, and $u_3=yy$ provide the partition of Lemma \ref{partition_qo}.
At the same time,  $v_0w_0=xxxxxyy\not\in O_3(S)$ while $w=xxxxyy$ is a $3$-overlap with $\rho_2(w)=xxxx$ and $\rho_1(w)=xxx\in S$. The paths $u_1=x$, $v_1=xx$, $u_2=e$, $v_2=x$, and $u_3=yy$ 
provide the partition of Lemma \ref{partition_o}. This example shows that it is possible that $(w_0,v_0)\in QO_n(S)$ while $v_0w_0\not\in O_n(S)$.
\end{ex}

Let us introduce the following notation. Given $n\in \mathbb{N}$, 
\[
\begin{aligned}
maxo_n(S)&=sup\{len(w)\mid w\in O_n(S)\},\\
mino_n(S)&=inf\{len(w)\mid  w\in O_n(S)\},\\
maxqo_n(S)&=sup\{len(w)\mid \exists v\mbox{ such that }(w,v)\in QO_n(S)\},\\
minqo_n(S)&=inf\{len(w)\mid \exists v\mbox{ such that }(w,v)\in QO_n(S)\}.
\end{aligned}
\]

By definition, we set $mino_n(S)=+\infty$ and $maxo_n(S)=-\infty$ if $O_n(S)$ is empty and $minqo_n(S)=+\infty$ and $maxqo_n(S)=-\infty$ if $QO_n(S)$ is empty.
Note that under this convention we have $mino_n(S)\ge n+1$ and $maxo_n(S)\le len(S)n-n+1$, where $len(S)$ denotes the maximal length of the paths in $S$.

Now we are going to prove a Theorem  that allows to estimate the values of $maxo_n(S)$ and $mino_n(S)$ using $maxqo_n(S)$ and $minqo_n(S)$.

\begin{theorem}\label{qo_o}  Given $n\ge 0$, for any reduced set $S$ we have
\begin{itemize}
\item $maxqo_n(S)\le maxo_n(S)-1$,
\item $minqo_n(S)\ge mino_n(S)-len(S)+1$.
\end{itemize}
\end{theorem}
\begin{Proof} The result is obvious for $n=0$. For $n\ge 1$, we are going to prove the following assertion. If $(w,v)$ is an $n$-quasioverlap, then there is $v'\in Q_{>0}$ such 
that $v'\mid_rv$ and $v'w$ is an $n$-overlap.
Since it follows easily from the definition of $1$-overlap that $len(v)<len(S)$, we have $$len(v'w)-len(S)+1\le len(w)+len(v)-len(S)+1\le len(w)\le len(v'w)-1$$ for such $v'$. 
Thus, after proving the existence of $v'$ we will be done.

More precisely, we will prove the following statements by induction on $n$. If $(w,v)\in QO_n(S)$, then there exists $v'\in Q_{>0}$ such that $v'\mid_r v$, $v'w\in O_n(S)$, $v'\rho_i^v(w)\mid_l \rho_i(v'w)$ for odd $i$, and $\rho_i(v'w)\mid_l v'\rho_i^v(w)$ for even $i$.

For $n=1$, we define $v'=v$.

Let us now consider the case $n=2$. If $vw$ cannot be presented in the form $vw=u'su$ with $u,u'\in Q_{>0}$ and $s\in S$, then we can take $v'=v$ and obtain the $2$-overlap $wv$ 
satisfying all the required conditions.
Suppose that it is possible to write $vw=u'su$ with $u,u'\in Q_{>0}$ and $s\in S$. Choose such a presentation with minimal $len(u)$. It follows from the 
definition of $2$-quasioverlap that $len(u')<len(v)$, i.e. there is $v'\in Q_{>0}$ such that $v=u'v'$.
It is easy to see that $v'w$ is a $2$-overlap with $\rho_1(v'w)=s$ satisfying all the required conditions.

Let us now prove the inductive step.  Suppose that the assertion above is true for all $(n-1)$-quasioverlaps. We will prove it for $n$ by induction on $len(v)$. 
Since the assertion is obvious for any $(w,v)\in QO_n(S)$ with $len(v)=1$, we may assume that, when we try to prove the assertion for some $n$-quasioverlap, we have already 
proved it for all $n$-quasioverlaps with length smaller than $len(v)$.

Let us consider $(w,v)\in QO_n(S)$ and denote by $w'$ the path $\rho_{n-1}^v(w)$. Since $(w',v)\in QO_{n-1}(S)$, we can apply the induction hypothesis. Thus, there exists $\tilde v\in Q_{>0}$ 
such that $\tilde v\mid_r v$, $\tilde vw'\in O_{n-1}(S)$, $\tilde v\rho_i^v(w')\mid_l \rho_i(\tilde vw')$ for odd $i$, and $\rho_i(\tilde vw')\mid_l \tilde v\rho_i^v(w')$ for even $i$. 
Note also that $\rho_i^v(w')=\rho_i^v(w)$ for $0\le i\le n-1$. Let us consider three cases:

1. $\rho_{n-2}(\tv w')\mid_S^a \tv w$. In this case we can simply define $v'=\tv$. It is easy to see that $v'w$ is an $n$-overlap with $\rho_{n-1}(v'w)=v'w'$ satisfying all the required conditions.

2. $\rho_{n-2}(\tv w')\mid_S \tv w$. In this case  $\tv\rho_{n-2}^v(w)\mid_l\rho_{n-2}(\tv w')$, $\tv\rho_{n-2}^v(w)\not=\rho_{n-2}(\tv w')$, and, hence, $2\nmid n$.

Suppose that $i$ is odd and $\rho_{i}(\tv w')\mid_S \tv\rho_{i+2}^v(w)$. Since 
\[\rho_i(\tv w')\mid_l\rho_{i+1}(\tv w'), \,\, \rho_{i+1}(\tv w')\mid_l\tv\rho_{i+1}^v(w'),\,\, \tv\rho_{i+1}^v(w') \mid_l \tv\rho_{i+2}^v(w),\text{ and }\rho_{i}(\tv w')\not=\rho_{i+1}(\tv w'),\] 
 we know that $\rho_{i}(\tv w')\not= \tv\rho_{i+2}^v(w)$. 

Also,
$\tv\rho_{i}^v(w)\mid_S^a \tv\rho_{i+2}^v(w) \text {  implies  } \tv\rho_{i}^v(w)\mid_S\rho_{i}(\tv w')\text{ and }\tv\rho_{i}^v(w)\not=\rho_{i}(\tv w').$
Thus, it follows from $\rho_{i-2}(\tv w')\mid_S^a \rho_{i}(\tv w')$ that $\rho_{i-2}(\tv w')\mid_S \tv\rho_{i}^v(w)$.

Now, the descending induction on $i$ gives us $\rho_{i}(\tv w')\mid_S \tv\rho_{i+2}^v(w)$ and $\tv\rho_{i}^v(w)\mid_S\rho_{i}(\tv w')$ for all odd $i$ such that $1\le i\le n-2$. 
As before, we have $\tv\rho_{1}^v(w)\not=\rho_{1}(\tv w')$. Consequently, $\tv\not=v$ and there is $v_0\in Q_{>0}$ such that $v=v_0\tv$.
Then $(\tv w,v_0)\in QO_{n-1}(S)$ and $len(v_0)<len(v)$. Thus we have $v_0'$ such that $v_0'\tv w\in O_n(S)$ satisfies all required conditions. It is easy to check that we can take $v'=v_0'\tv$ 
in this case.

3. There is a presentation $\tv w=\rho_{n-2}(\tv w')u'su$ with $u\in Q_{>0}$, $u'\in Q$ and $s\in S$. Let us choose such a presentation with minimal $len(u')$. 
Since $\rho_{n-2}(\tv w')\mid_l\tv\rho_{n-2}^v(w)$ and $\rho_{n-2}(\tv w')\not=\tv\rho_{n-2}^v(w)$, $n$ must be even.

It follows from the minimality of $u'$ that $\rho_{n-2}(\tv w')\mid_S^a \rho_{n-2}(\tv w')u's$.
Since $\rho_{n-2}(\tv w')\mid_S\tv w'$, we know that $\tv w'\mid_l\rho_{n-2}(\tv w')u's$ and, hence, $\rho_{n-2}(\tv w')u's\mid_S\tv w$.
Since $\tv\rho_{n-2}^v(w)\mid_S^a\tv w$, we have also $\tv\rho_{n-2}^v(w)\mid_S\rho_{n-2}(\tv w')u's.$

For even $2\le i\le n-2$,  we have $\tv\rho_{i-1}^v(w)\mid_l\rho_{i-1}(\tv w')\mid_l\rho_{i}(\tv w')\mid_l\tv\rho_{i}^v(w)$ and $\tv\rho_{n-1}^v(w)\mid_S\tv\rho_{i}^v(w)$. 
Thus, $\rho_{i}(\tv w')\mid_S\tv\rho_{i}^v(w)$ for such $i$.
Suppose that $\rho_{i}(\tv w')\mid_S\tv\rho_{i}^v(w)$ and $\rho_{i}(\tv w')\not=\tv\rho_{i}^v(w)$ for some $2\le i\le n-2$. Since  $\tv\rho_{i-2}^v(w)\mid_S^a\tv\rho_{i}^v(w)$, 
we have $\tv\rho_{i-2}^v(w)\mid_S\rho_i(\tv w')$.
Moreover, $\rho_{i-2}(\tv w')\mid_S^a\rho_i(\tv w')$ implies $\rho_{i-2}(\tv w')\not=\tv\rho_{i-2}^v(w)$.

Thus, the descending induction on $i$ gives us  $\rho_i(\tv w')\mid_S\tv\rho_i^v(w)$ and $\tv\rho_{i-2}^v(w)\mid_S\rho_i(\tv w')$ for even $i$, $2\le i\le n-2$.
In particular, we get $\tv\mid_S\rho_2(\tv w')$. On the other hand, there is $r\in S$ such that $r\mid_r\rho_2(\tv w')$ and $r\not=\rho_2(\tv w')$. 
As a consequence,$\rho_2(\tv w')=v_0'r$ and $\tv=v_0'v_0$ for some $v_0,v_0'\in Q_{>0}$. In particular, $len(v_0)<len(v)$.
It follows from our arguments that $(w,v_0)\in QO_n(S)$ with $v_0'r_{n-1}^{v_0}(w)=\rho_{n-2}(\tv w')u's$. We thus have got $v'$ such that $v'w\in O_n(S)$ satisfies all required conditions.
\end{Proof}

\begin{coro}\label{sum_o} For any reduced set of paths $S$, 
\[ maxo_{n+m}(S)\le maxo_n(S)+maxo_m(S)-1\] 
and 
\[mino_{n+m}(S)\ge mino_n(S)+mino_m(S)-len(S)+1.\]
\end{coro}
\begin{Proof} The proof follows from 
Theorem  \ref{qo_o} and the fact that any $w\in O_{n+m}(S)$ can be represented in the form $w=w'w''$ with $w'\in O_n(S)$ and $w''\in QO_m(S)$.
\end{Proof}

From now on we fix an admissible order $\ge$ on the set of paths in $Q$, \cite{gff}.
More precisely, this means that there is a well order $\ge$ such that for any paths $p,q,u,v \in Q$, 
\begin{itemize}
 \item if $p\ge q$, then $upv\ge uqv$ if the products are paths.
 \item $p\ge q$ if $q\mid p$.
\end{itemize}

Given a linear space $V$, its basis $B$, and $x\in V$, we call the sum $\sum\limits_{i=1}^m\alpha_ib_i$ a {\it reduced expression} of $x$ as a linear combination of the elements of $B$ if $\alpha_i\in\kk^*$, $b_i\in B$ for $1\le i\le m$, $b_i\not=b_j$ for $1\le i<j\le m$, and $x=\sum\limits_{i=1}^m\alpha_ib_i$.

For $x\in \kk Q$,  $tip(x)$ is maximal path of $Q$, with respect to the order $\ge$,
appearing in the reduced expression of $x$ as a linear combination of paths.
A subset $G\subset I$ is called a {\it Gr\"obner basis} of $I$ if for any $x\in I$, there exists $g\in G$ such that 
$tip(g)\mid tip(x)$.

From now on we fix $A=\kk Q/I$, where $I$ is an ideal of $\kk Q$ that has a finite Gr\"obner basis $G$.

We will use the notation of Green and Solberg in {\cite{gs}}.
Let $X$ be a graded $A$-module and let $P_0(X)\xrightarrow{\mu_X}X$ be its minimal graded projective cover. Suppose also that $X$ is finitely presented. 
The projective module $P_0(X)$ can be presented in the form $P_0(X)=\oplus_{i\in T_0}f_i^0A$, where $T_0$ is a finite set and, for any $i\in T_0$ 
there exist $e_i\in Q_0$, $m_i\in\mathbb{Z}$, and an isomorphism of graded modules  $f_i^0A\cong e_iA[m_i]$ that sends $f_i^0$ to $e_i$.

Consider the graded space $V_X=\oplus_{i\in T_0}f_i^0\kk Q$. The
set whose elements  are of the form $f_i^0p$, where $i\in T_0$ and $p$ is a path ending in $e_i$, is a basis of $V_X$, we denote this set by $B_X$ and introduce the 
following well order on it. We set $f_i^0p\ge f_j^0q$ if either
$p>q$ or $p=q$ and $i\ge j$. We say that $f_i^0p$ {\it divides} $f_j^0q$ {\it on the right} if $i=j$ and $p\mid_rq$. For
$x\in V_X$, we write $tip(x)$ for the maximal element of $B_X$, with respect to the order $\ge$, appearing in the reduced
expression of $x$ as a linear combination of elements of $B_X$. The set $x_1,\dots,x_l$ of nonzero elements of $V_X$ is
called {\it right tip reduced} if $tip(x_i)$ does not divide $tip(x_j)$ on the right  for any $1\le i,j\le l$, $i\not=j$.

By \cite[Proposition 5.1]{gs}, there are finite sets $T_1$ and $T_1'$, elements $h_i^1\in V_X$ ($i\in T_1$), and elements ${h_i^1}'\in V_X$ ($i\in T_1'$) such that
\begin{enumerate}
\item any element $h$ in the kernel of the composition $V_X\twoheadrightarrow P_0(X)\xrightarrow{\mu_X}X$ can be uniquely represented in the form $h=\sum_{i\in T_1}f_i+\sum_{i\in T_1'}f_i'$, 
where $f_i\in h_i^1\kk Q$ and $f_i'\in {h_i^1}'\kk Q$,
\item for any element $h\in\{h_i^1\}_{i\in T_1}\cup \{{h_i^1}'\}_{i\in T_1'}$ there exists $e_h\in Q_0$ such that $he_h=h$,
\item ${h_i^1}'\in \oplus_{j\in T_1}h_j^0I$ for any $i\in T_1'$,
\item the set $\{h_i^1\}_{i\in T_1}\cup \{{h_i^1}'\}_{i\in T_1'}$ is right tip reduced.
\end{enumerate}
Moreover, it is clear that all the elements in the set $\{h_i^1\}_{i\in T_1}\cup \{{h_i^1}'\}_{i\in T_1'}$ can be chosen homogeneous. 
Let us define $\bar P_1(X)=\oplus_{i\in T_1}f_i^1A$, where, for any $i\in T_1$ there exists an isomorphism of graded modules $f_i^1A\cong e_{h_i^1}A[-deg(h_i^1)]$ 
that sends $f_i^1$ to $e_{h_i^1}$. Here $deg(h_i^1)$ denotes the degree of $h_i^1$. Note that the order in the basis of $V_X$ induces an order on the set $\{f_i^1\}_{i\in T_1}$. 
In this way we obtain a graded projective presentation
$$
\bar P_1(X)\xrightarrow{\bar d_0(X)}P_0(X)\xrightarrow{\mu_X}X
$$
of $X$, where $\bar d_0(X)$ sends $f_i^1$ to the class of $h_i^1$ in $P_0(X)$ for any $i\in T_1$.
Let also
\[
\dots \xrightarrow{d_n(X)} P_n(X)\xrightarrow{d_{n-1}(X)}\dots \xrightarrow{d_0(X)}P_0(X)(\xrightarrow{\mu_X}X)
\]
be the minimal graded $A$-projective resolution of $X$. 

 As before,  given a  graded $A$-module $M$, we will denote by $M_j$ its $j$-th homogeneous component. The next Theorem will be crucial for us. 
It allows to estimate the degrees of the generators of the terms of the resolution constructed using the algorithm from \cite{gs}.

\begin{theorem}\label{qo}
Let us fix the notation as above. If $k$ and $l$ are the minimal and the maximal degrees of $f_i^1$, $i\in T_1$, then, for any $n\ge 1$,
$P_{n}(X)=\left(\bigoplus\limits_{j=k+minqo_n(tip(G))}^{l+maxqo_n(tip(G))}P_{n}(X)_j\right)A$.
\end{theorem}
\begin{Proof} The algorithm described in \cite[Section 3]{gs} gives a graded projective resolution
\[
\dots \xrightarrow{\bar d_n(X)} \bar P_n(X)\xrightarrow{\bar d_{n-1}(X)}\dots \xrightarrow{\bar d_0(X)}P_0(X)(\xrightarrow{\mu_X}X),
\]
of $X$. Due to this algorithm $\bar P_n(X)$ can be presented in the form $\bar P_n(X)=\oplus_{i\in T_n}f_i^nA$ and $\bar d_{n-1}(X)$ sends 
$f_i^n$ to $h_i^n\in \oplus_{i\in T_{n-1}}f_i^{n-1}A$, where all these elements satisfy the following property.
If $tip(h_i^n)=f_j^{n-1}p$, then $tip(h_j^{n-1})p\in Q_S$ and if $q\mid_ltip(h_j^{n-1})p$, $q\not=tip(h_j^{n-1})p$, then $q$ is not $S$-path.

For $i\in T_n$, let us introduce $p_i\in Q$ and $t(i)\in T_{n-1}$ by the equality $tip(h_i^n)=f_{t(i)}^{n-1}p_i$. 
It is possible to prove inductively on $n$ that $(w,v)\in QO_n(S)$ for $w=p_{t^{n-2}(i)}\dots p_{t(i)}p_i$ and some $v\in Q_{>0}$ such that 
$v\mid_rtip\left(h_{t^{n-1}(i)}^1\right)$.
Consequently, the required statement follows from the equality
\[
deg\left(f_i^n\right)=deg\left(f_{t^{n-1}(i)}^1\right)+\sum\limits_{m=0}^{n-2}len\left(p_{t^m(i)}\right)=deg\left(f_{t^{n-1}(i)}^1\right)+len(w).
\]
\end{Proof}

In other words, for any $n\ge 1$, we are able to deduce the possible degrees of the generators of $P_{n}(X)$ from those of the generators of $\bar P_1(X)$ and the 
lengths of $n$-quasioverlaps for $tip(G)$.

We obtain two corollaries.

\begin{coro}\label{over}
Let $A=\kk Q/I$, $X$ be a graded finitely presented $A$-module and let $G$ be a homogeneous Gr\"obner basis for $I$ such that $tip(G)$ is reduced. Let $P_{\bullet}(X)$  be a minimal graded $A$-projective resolution of the 
graded $A$-module $X$ and $\bar P_1(X)$ be as above. If $k$ and $l$ are the minimal and the maximal degrees of $f_i^1$, $i\in T_1$, then, for any $n\ge 1$,
$P_{n}(X)=\left(\bigoplus\limits_{j=k+mino_n(tip(G))-len(tip(G))+1}^{l+maxo_n(tip(G))-1}P_{n}(X)_j\right)A$.
\end{coro}
\begin{Proof} It follows directly from Theorems  \ref{qo_o} and \ref{qo}.
\end{Proof}

\begin{coro} Let $A=\kk Q/I$ where $I$ has a homogeneous Gr\"obner basis $G$ such that $len(tip(G))\le s$ and $maxo_2(tip(G))\le s+1$. 
If additionally $mino_1(tip(G))=s$, then the algebra $A$ is $s$-Koszul.
\end{coro}
\begin{Proof} It follows from Corollaries \ref{sum_o} and \ref{over} since $\bar P_1(A_0)=\oplus_{\alpha\in Q_1}e_{\alpha}A[-1]$, where $e_{\alpha}$ is the starting vertex of the arrow $\alpha$.
More precisely, we get by induction on $i\ge 3$ that
$$maxo_i(tip(G))\le maxo_{i-2}(tip(G))+maxo_2(tip(G))-1\le \chi_s(i-2)+s=\chi_s(i).$$ 
If additionally $mino_1(tip(G))=s$, then all the elements of $tip(G)$ have length $s$. Then we get that $P_2(A_0)$ is generated in degree $s+1$ and that if the minimal generating degree for $P_i(A_0)$ is $m$, then the minimal generating degree for $P_{i+2}(A_0)$ is not less than $m+s$.
Then we get by induction that $P_i(A_0)$ is generated in degree $\chi_s(i)$.
\end{Proof}

\vspace{0.5 cm} 

\noindent\textsf{Eduardo N. Marcos: IME-USP (Departamento de Matem\'atica),
  Rua Mat\~ao 1010 Cid. Univ., S\~ao Paulo, 055080-090, Brasil.}
	
\noindent\emph{enmarcos@ime.usp.br}\\

\noindent\textsf{Andrea Solotar: IMAS and Dto de Matem\'{a}tica, Facultad de Ciencias Exactas y Naturales,
Universidad de Buenos Aires, Ciudad Universitaria, Pabell\'{o}n 1,
(1428) Buenos Aires, Argentina.}

\noindent\emph{asolotar@dm.uba.ar}\\

\noindent\textsf{Yury Volkov: Saint-Petersburg State University, Universitetskaya nab. 7-9, St. Petersburg, Russia.\\
Dto de Matem\'atica, Instituto de Matem\'atica e Estat\'istica, Universidade S\~ao Paulo, Rua de Mat\~ao 1010, Cidade Universit\'aria, S\~ao Paulo-SP, 055080-090, Brasil.}

\noindent\emph{ wolf86\_666@list.ru}

\end{document}